\documentclass[11pt]{article}
\usepackage{amsmath,amsthm,amssymb}
\usepackage{graphicx}
\usepackage{mathrsfs}
\usepackage{hyperref}
\usepackage{mathtools}
\usepackage{enumitem}
\hypersetup{colorlinks=true,linkcolor=blue,filecolor=blue,citecolor=blue}
\usepackage{cleveref}
\usepackage{cite}
\usepackage{makecell}
\usepackage{array}
\usepackage[defaultlines=3,all]{nowidow}
\usepackage{microtype}
\usepackage{multirow}
\usepackage{mdwtab}
\usepackage{comment}
\newtheorem{theorem}{Theorem}[section]
\newtheorem{lemma}[theorem]{Lemma}

\newtheorem{definition}[theorem]{Definition}

\newtheorem{claim}{Claim}
\newtheorem{fact}{Fact}
\newtheorem{observation}{Observation}
\newtheorem{question}{Question}
\newtheorem{subquestion}{Subquestion}

 \usepackage{todonotes} 

\usepackage{thmtools}
\usepackage{cleveref}
\crefname{subquestion}{Subquestion}{Subquestions} 
\Crefname{subquestion}{Subquestion}{Subquestions} 

\usepackage{pifont}      

\usepackage{algorithm}  
\usepackage{algpseudocode}  
\usepackage{amsmath}  

\usepackage{float}  

\begin{document}

\title{\bf Computational results on semistrong edge coloring of graphs
}
\author{Yuquan Lin 
and Wensong Lin
\footnote{Corresponding author. E-mail address: wslin@seu.edu.cn}\\
{\small School of Mathematics, Southeast University, Nanjing 210096, P.R. China}}
\date{}
\maketitle

\vspace*{-1cm} 
\medskip

\begin{abstract}
\medskip
The semistrong edge coloring, as a relaxation of the well-known strong edge coloring, can be used to model efficient communication scheduling in wireless networks.
An edge coloring of a graph $G$ is called \emph{semistrong} if every color class $M$ is a matching such that every edge of $M$ is incident with a vertex of  degree 1 in the subgraph of $G$ induced by the endvertices of edges in $M$.
The  \emph{semistrong  chromatic index} $\chi_{ss}'(G)$ of  $G$ is the minimum number of colors required for a semistrong edge coloring.
In this paper, we 
prove that the problem of	determining whether a graph $G$ has a semistrong edge coloring with $k$ colors is polynomial-time solvable for   $k\le2$ and
 is NP-complete for  $k\ge3$.
For trees, we develop  a  polynomial-time algorithm  to determine the  semistrong chromatic index exactly.
\end{abstract}


\section{Introduction}
\label{sec:intro}
The concept of a $\mathscr{P}$-matching
was introduced in \cite{GHHL2005}  to unify and extend various types of matchings studied in graph theory.
Let $G$ be a  graph, a \emph{matching} $M$ of  $G$ is  a  set of pairwise disjoint edges, and  the subgraph  of $G$
 induced by the vertices covered by $M$ is denoted by
$G_M$.
Given a graph property $\mathscr{P}$, a matching $M$ of $G$ is called a
\emph{$\mathscr{P}$-matching} if  $G_M$ satisfies the property $\mathscr{P}$.
From the   strongest to weakest property $\mathscr P$, 
specific examples include (See \Cref{fig:matching}):
\begin{itemize}
    \item \emph{induced matching} (also known as \emph{strong matching}): every vertex of $G_M$ is pendant, where a \emph{pendant vertex} is a vertex of degree $1$.
     \item \emph{semistrong matching}: every edge of $M$ is pendant in $G_M$, where a \emph{pendent edge} is an edge incident to a pendent  vertex.     
    \item \emph{uniquely restricted matching}: $G_M$ has a unique perfect matching  (Note that $M$  itself forms a perfect matching of $G_M$).
\end{itemize}
\begin{figure}[H]
	\centering
	\resizebox{16cm}{2.6cm}{\includegraphics{matching.pdf}}
    \caption{Different types of matching (heavy edges), from weaker to stronger.}
    \label{fig:matching}
\end{figure}

\begin{figure}[h!t]
	\begin{center}	\includegraphics[width=.8\textwidth]{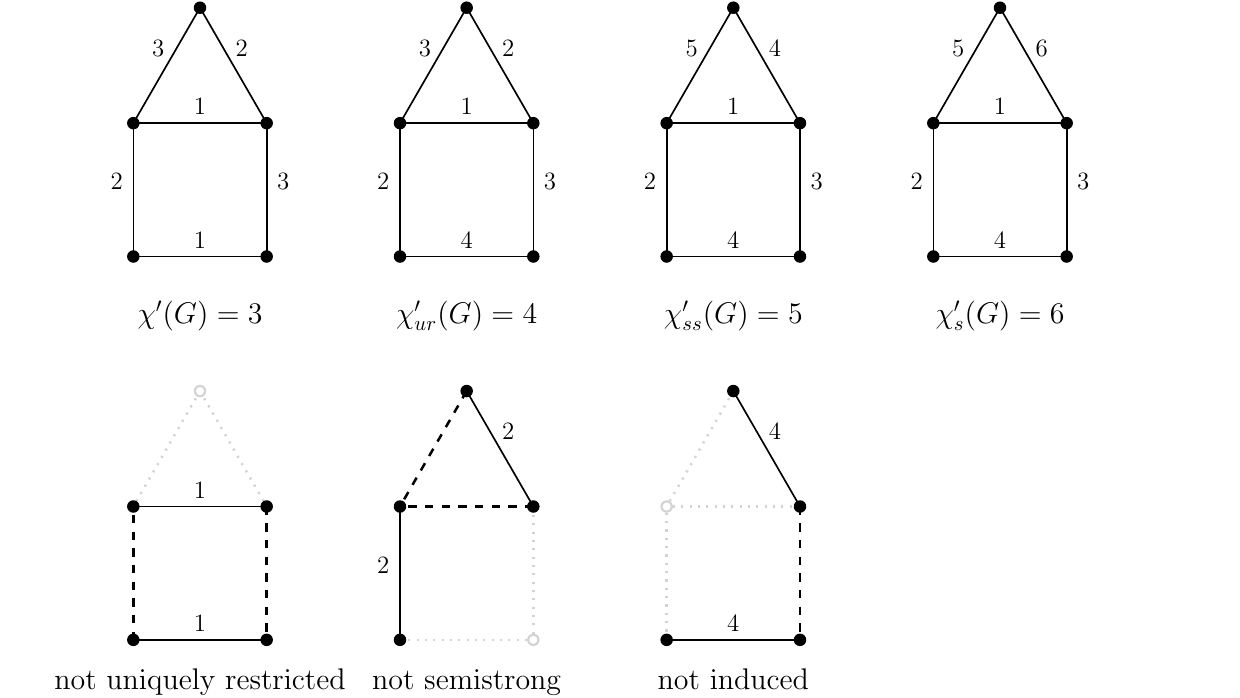}
		\caption{Different types of edge coloring (proper, uniquely restricted, semistrong, strong).}
		\label{fig:house}
	\end{center}
\end{figure}

Each type of matching naturally leads to a corresponding notion of edge coloring.
Specifically, 
the  \emph{proper edge coloring},
\emph{uniquely restricted edge coloring} \cite{BRS2017},  \emph{semistrong edge coloring} \cite{GH2005}, and  \emph{strong edge coloring} \cite{FJ1983} of a graph $G$ is an assignment of colors to its edges such that each color class is a matching,  a uniquely restricted matching, a semistrong matching, and an induced matching, respectively (See \Cref{fig:house}). The \emph{chromatic index} $\chi'(G)$, 
\emph{uniquely restricted  chromatic index} $\chi_{ur}'(G)$,  \emph{semistrong  chromatic index} $\chi_{ss}'(G)$, and  \emph{strong  chromatic index}  $\chi_{s}'(G)$ of a graph $G$ are defined as the minimum number of colors required for the corresponding colorings.

Recall that a matching $M$ is uniquely restricted if and only if there is no $M$‐alternating cycle (a cycle that alternates between edges in $M$ and edges not in $M$ ) in $G$ \cite{GHL2001}. 
The uniquely restricted edge coloring is a strengthening of the   \emph{acyclic edge coloring} introduced in \cite{F1978}, which is an edge coloring of a graph $G$ such that each cycle of $G$ gets at least 3 colors.
The \emph{acyclic  chromatic index} $a'(G)$ is defined as the minimum number of colors required in an acyclic edge coloring of $G$.
Due to the  strength of the constraints corresponding to these colorings, every graph $G$ satisfies
\begin{equation}\label{ineq}
    \chi'(G)\le a'(G)\le \chi_{ur}'(G)\le \chi_{ss}'(G)\le \chi_{s}'(G).
\end{equation}

The well-known Vizing's theorem shows that $\chi'(G)$ is either $\Delta$ or $\Delta+1$  for any simple graph $G$ with maximum degree $\Delta$ \cite{V1964,G1967}.
 Nevertheless, Holyer \cite{H1981} proved that  the determination
 of $\chi'(G)$ is NP-hard.
 Regarding acyclic edge coloring, 
 Alon and Zaks \cite{AZ2002} established its NP-completeness by showing that  deciding whether $ a'(G)\le3$ holds for an arbitrary graph $G$ is NP-complete.
 Strong edge coloring, being a more restrictive and special case, is intuitively harder than the previous coloring problems.
  Mahdian \cite{M2002} proved that  deciding whether a graph $G$ admits a strong edge coloring using $k$ colors  is polynomial-time solvable for   any  integer $k\le3$, but 
  is NP-complete for any  integer $k\ge4$.
  
Motivated by two famous open conjectures on $a'(G)$  and $\chi_{s}'(G)$ (Fiamčik\cite{F1978} and Alon et al \cite{ASZ2001} conjectured $a'(G)\le \Delta+2$ holds	for any graph $G$;  Erdős and Nešetřil \cite{EN1989} conjectured $\chi_{s}'(G)\le 1.25 \Delta^2$ holds	for any graph $G$),  
the chain of inequalities in~(\ref{ineq})
naturally leads to further investigation of $\chi_{ur}'(G)$ and  $\chi_{ss}'(G)$, which are sandwiched between the acyclic and strong chromatic indices.

This paper focus on  the computational complexity and algorithm on  semistrong edge coloring.
Before presenting specific results,   we would like to point out that the semistrong edge coloring problem serves as a useful model for efficient communication scheduling in wireless networks.

\medskip
\noindent{\textbf{Application.}} 
Wireless devices, such as sensors and controllers, are extensively deployed in urban sensing and smart building environments to periodically synchronize status information or exchange local data, hereafter referred to as \emph{communication}. 
This network can be modeled as a graph 
$G=(V,E)$, where vertices represent wireless devices and edges denote pairs of   nearby  devices that  require communication.
To mitigate interference, each device is constrained to participate in at most one communication at any given time. 
Consequently, the set of communication tasks scheduled concurrently corresponds to a matching in the graph $G$.
In order to conserve energy, for each communicating pair, only one device acts as the initiator (called the \emph{caller}), responsible for sending connection requests and occupying channel resources, while the other device (called the \emph{callee})  solely receives. Notably, after communication is established, data exchange is bidirectional, but the caller bears most of the resource consumption.

Due to the broadcast nature of wireless signals, when a caller initiates a communication request, the signal can be received by all nearby devices, not just the intended callee. 
Limiting each callee to receive at most one incoming request at any given time helps reduce physical signal interference and ensures reliable connection establishment.
This spatial constraint leads to the communication tasks within each time slot forming a semistrong matching in the graph $G$, where the pendent vertex of each edge in this matching corresponds to the callee.
Under these conditions, the minimum time required to complete all communication tasks is characterized by the semistrong chromatic index of $G$.



\subsection*{Our results}

   In this paper, we first establish the  computational complexity of semistrong edge coloring problem by proving
    the following.
  \begin{theorem}\label{Main-th-NPC}
  	For any fixed integer $k$, the problem of	determining whether a graph $G$ has a  semistrong edge coloring using at most $k$ colors
    is polynomial-time solvable for  $k\le2$ and
  	is NP-complete for  $k\ge3$.
  \end{theorem}

Then, we address the semistrong edge coloring problem for trees. 
It was pointed out in \cite[Corollary 2]{LMS2024} that any tree $T$ with  maximum degree $\Delta$ has semistrong chromatic index either $\Delta$ or $\Delta+1$, and it is easy to get a semistrong edge coloring of $T$ using $\Delta+1$ colors using the algorithm provided in \cite[Algorithm 2]{HL2017}.
Based on this fact, we develop    a dynamic programming algorithm to determine whether a tree $T$ admits a  semistrong edge coloring with $\Delta$ colors. Our result is as follows.
 

  \begin{theorem}\label{Main-th-tree}
	Let $T$ be a tree with maximum degree $\Delta$.
	Then one can determine whether the semistrong chromatic index of $T$ is equal to either $\Delta$ or $\Delta+1$  in $O(\Delta^{6}n)$ time.
\end{theorem}

 Remark that if the semistrong chromatic index of $T$  is $\Delta$, 
 a semistrong edge coloring of $T$ using $\Delta$ colors can be readily derived based on our algorithm (See 	\Cref{alg:main}).
\bigskip

The remainder of this paper is organized as follows.
In the next section,
 we start by introducing some notation and preliminaries. 
 Then, \Cref{Main-th-NPC} and \Cref{Main-th-tree} are proved in \Cref{sec:NPC} and \Cref{sec:tree}, respectively.
 Finally, we discuss some further research directions in \Cref{sec:summary}.


 \section{Notation and preliminaries}

This paper considers only  finite undirected simple graphs.
Given a graph $G=(V(G),E(G))$.
 For $v\in V(G)$, let $N_{G}(v)=\{u\in V(G): uv\in E(G)\}$ denote the open neighborhood of $v$ and  $d_{G}(v)=|N_{G}(v)|$ be the degree of $v$. Let  $\Delta(G)=\max_{v\in V(G)}d_{G}(v)$ denote the maximum degree of $G$. 
We call $G$ a {\em $k$-regular graph} if
 $d_{G}(v)=k$ for every $v\in V(G)$.
 For any $X\subseteq V(G)$, we use $G-X$ to denote the subgraph of $G$ obtained by deleting vertices in $X$. 
  And for any $E'\subseteq E(G)$, we use $G\setminus E'$ to denote the subgraph of $G$ obtained by deleting edges in $E'$.
  In particular, we use the abbreviations  $G-v$  and  $G\setminus e$ for the graphs $G-\{v\}$ and $G\setminus \{e\}$, respectively.

 
  The  distance between two edges $e$ and $f$ in $G$ is equal to the distance between their corresponding vertices in the line graph $L(G)$ of $G$, denoted by $d_{G}(e,f)$.  For any two edges $e,f\in E(G)$,
  we say that  $f$ is a {\em $1$-neighbor} (resp.  {\em $2$-neighbor}) of $e$ if  $f$ and $e$ are at distance 1 (resp.  2). 

   Let $\phi$ be an edge coloring of $G$. For any edge subset $E'\subseteq E(G)$, we denote by $\phi(E')$ the set of colors that appear on the edges in $E'$.
For an edge $e=uv$ in $G$, we say that {\em $u$ is a 1-vertex of $e$ under $\phi$} (or say, {\em $e$ has a 1-vertex $u$ under $\phi$}) if $u$ is a vertex of degree 1 (or say, a pendant vertex) in the subgraph of $G$ induced by the endvertices of edges in $G$ that is colored with $\phi(e)$. It is obvious that $\phi$  is a semistrong edge coloring of $G$  if and only if every edge in $G$ has at least one 1-vertex  under $\phi$. 
 Given a positive integer $k$,  if  $G$ has a  (semistrong) edge coloring $\phi$ using $k$ colors, then we say that   $G$ is  {\em (semistrongly) $k$-edge-colorable} and  $\phi$ is a \emph{(semistrong) $k$-edge-coloring}  of  $G$.
 

	\section{NP-completeness}\label{sec:NPC}

	In this section, we focus on the computational complexity of semistrong edge coloring.
	The decision version of the semistrong edge coloring problem is stated as follows:
    
\medskip	
	\noindent{\textbf{Semistrong Edge Coloring Problem}}
\\
	{\em Instance}: A graph $G$ and a positive integer $k$.\\
	{\em Question}: Does $G$ admit a semistrong $k$-edge-coloring?

\medskip		
For $k\le2$, the above problem is trivial as the only semistrongly $2$-edge-colorable graphs are disjoint unions of paths of order at most 5.
In the following, we will show  this problem is NP-complete  for any other values of $k$ by giving a polynomial reduction from the following edge coloring problem,  which is known to be NP-complete for every fixed interger $k\ge3$ \cite{LG1983}:

\medskip
\noindent{\textbf{Edge Coloring Problem of Regular graphs}}
\\
{\em Instance}: A $k$-regular graph $G$ with $k\ge3$.\\
{\em Question}: Does $G$ admit a $k$-edge-coloring?

\medskip

Now, we prove the the following theorem to complete the proof of \Cref{Main-th-NPC}.
	
\begin{theorem}\label{Main-th}
For any fixed integer $k\ge3$,	Semistrong Edge Coloring Problem is NP-complete. 
\end{theorem}
\begin{proof}
	It is clear that the problem is in NP. 
For each	 $k\ge3$ and each $k$-regular graph $G$, we will construct another graph $H$ with maximum degree $k$ such that $G$ has a $k$-edge-coloring if and only if   $H$ has a semistrong $k$-edge-coloring.
According to the value of $k$, we divide the proof into the following two cases. 	

\textbf{Case 1}. $k\neq4$.

Let  $G$ be a $k$-regular graph  with $k\neq4$. We  construct a graph $H$  as follows: if $k\ge3$ is odd, we  replace  each edge $uv\in E(G)$ with the graph $B_{k}[u,v]$ shown in Figure  \ref{fig:odd}; and if  $k\ge6$ is even, we  replace  each edge $uv\in E(G)$ with the graph $Q_{k}[u,v]$ shown in Figure  \ref{fig:even}. 
It is obvious that $H$  is a graph with maximum degree $k$.
Denote by  $B_{k}(u,v)$ the graph  $B_{k}[u,v]-\{u,v\}$ and  $Q_{k}(u,v)$ the graph  $Q_{k}[u,v]-\{u,v\}$, respectively.
 See Figures \ref{fig:odd} and \ref{fig:even} for the names of vertices in $B_{k}[u,v]$ and $Q_{k}[u,v]$, respectively.

 \begin{figure}[htbp] 
	\centering	
	\begin{minipage}[t]{7cm}
		\centering
		\resizebox{5.6cm}{3.6cm}{\includegraphics{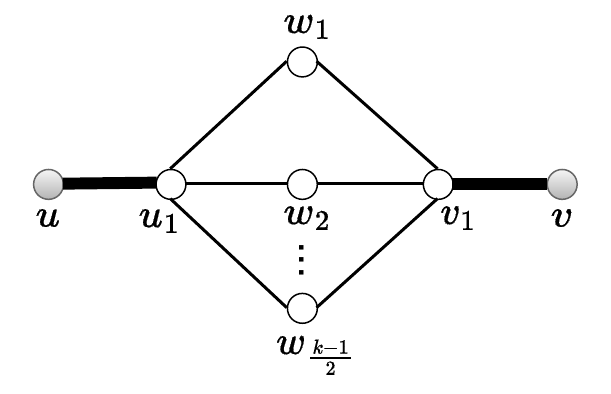}}
	\caption{The graph $B_{k}[u,v]$.}
	\label{fig:odd}
	\end{minipage}
	\begin{minipage}[t]{7cm}
		\centering
		\resizebox{6.8cm}{4.3cm}{\includegraphics{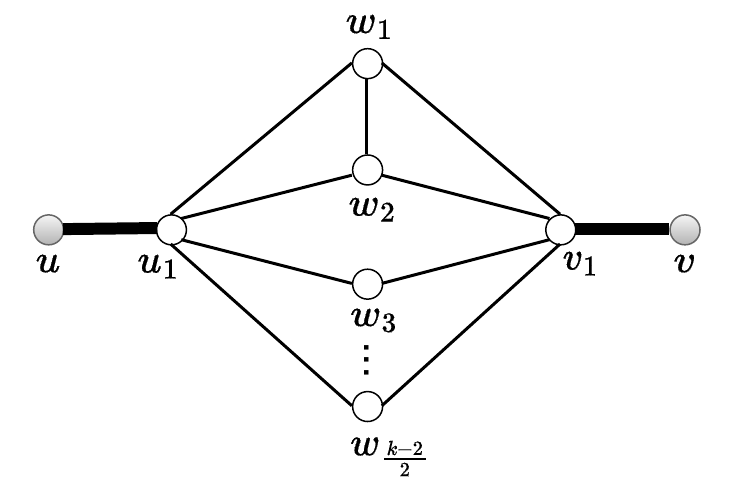}}
	\caption{The graph $Q_{k}[u,v]$.}
	\label{fig:even}
	\end{minipage}
\end{figure}

When $k$ is odd,   each vertex $w_{i}$ ($1\le i\le \frac{k-1}{2}$) is of degree 2 in  the graph $B_{k}[u,v]$ and so $B_{k}[u,v]$
 has exactly $k+1$ edges. Any  two distinct edges in $B_{k}(u,v)$ 
 must receive different colors in any semistrong edge coloring of $B_{k}[u,v]$ as they 
  are  at distince one or in a  common $4$-cycle.

When   $k$ is even,  both $w_{1}$ and $w_{2}$ are of degree 3 in   $Q_{k}[u,v]$  and  if  $k\ge 8$ then $w_{i}$ is of degree 2 in $R[u,v]$ for any $3\le i\le \frac{k-2}{2}$, 
 the graph  $Q_{k}[u,v]$ also has  $k+1$ edges.  
 Under any semistrong edge coloring of  $Q_{k}[u,v]$, the color of the edge $w_{1}w_{2}$ must be different  from all the other edges in $Q_{k}[u,v]$ and
 any  two distinct edges in $Q_{k}(u,v)$ 
 must receive different colors.

Let $\phi$ be a   $k$-edge-coloring of $G$.
For each $uv\in E(G)$, we first color the two edges $uu_{1}$ and $vv_{1}$ with the same color $\phi(uv)$ and then use the remaining $k-1$ colors to
color the other $k-1$  edges  in $B_{k}[u,v]$ (resp.  $Q_{k}[u,v]$).
	This results in an edge coloring $\psi$ of $H$. We next prove that $\psi$ is a semistrong $k$-edge-coloring of $H$.
Let $e$ be an edge in $H$.  
Without loss of generality, we may assume that 
$e$ belongs to the subgraph $B_{k}[u,v]$  (resp.  $Q_{k}[u,v]$) of $H$ corresponding to the edge $uv\in E(G)$.
If $e$ is incident with  $u$ (resp. $v$), then $u_{1}$ (resp. $v_{1}$) is a 1-vertex of  $e$ under $\psi$. 
And if  $e$ is incident with  $w_{i}$ for some  $1\le i \le \frac{k-1}{2}$ (resp.  $1\le i \le \frac{k-2}{2}$), then   $w_{i}$ is a 1-vertex of  $e$ under $\psi$.  
Therefore,  $\psi$ is indeed a semistrong $k$-edge-coloring of $H$.

Conversely, let  $\psi$ be a semistrong $k$-edge-coloring of $H$. 
We first prove that for each $uv\in E(G)$, it holds that $\psi(uu_{1})=\psi(vv_{1})$ in its corresponding subgraph $B_{k}[u,v]$ (resp.  $Q_{k}[u,v]$) of $H$.
If not,  we may assume that
$\psi(uu_{1})=1$ and $\psi(vv_{1})=2$.
Recall that any  two distinct edges in $B_{k}(u,v)$  (resp.  $Q_{k}(u,v)$)
must receive different colors under $\psi$ and  the graph $B_{k}(u,v)$ (resp.  $Q_{k}(u,v)$) has exactly $k-1$ edges, we have $|\psi(B_{k}(u,v))|= k-1$ (resp.  $|\psi(Q_{k}(u,v)|= k-1$). 
It follows that $\{1,2\}\cap \psi(B_{k}(u,v))\neq\emptyset$ (resp.  $\{1,2\}\cap \psi(B_{k}(u,v))\neq\emptyset$ ). 
We may assume that 
$\psi(u_{1}w_{1})=\psi(vv_{1})=2$.  
Since $u$ is of degree  $k$ in $H$, there is exactly one edge $f\neq uu_{1}$ incident with $u$ such that $\psi(f)=2$.
Because  both  $f$ and $vv_{1}$ are at  distance two from the edge $u_{1}w_{1}$, the edge $u_{1}w_{1}$ has no $1$-vertex under $\psi$.
This is a contradition to the semistrong $k$-edge-coloring $\psi$.
Therefore,  we can obtain a  $k$-edge-coloring  $\phi$ of $G$ by setting $\phi(uv)=\psi(uu_{1})=\psi(vv_{1})$
 for each $uv\in E(G)$. 

\textbf{Case 2}. $k=4$.

Let $G$ be a $4$-regular graph.  We  construct a  graph $H$ with maximum degree $4$ by replacing each edge $uv\in E(G)$ by a graph $R[u,v]$, see Figure \ref{fig:4} for  the names of vertices and edges in $R[u,v]$.
 Next, we prove that  $G$ has a $4$-edge-coloring if and only if   $H$ has a semistrong $4$-edge-coloring.

  \begin{figure}[htbp]
 	\centering	
 	\begin{minipage}[t]{7.9cm}
 		\centering
 		\resizebox{7.9cm}{6.1cm}{\includegraphics{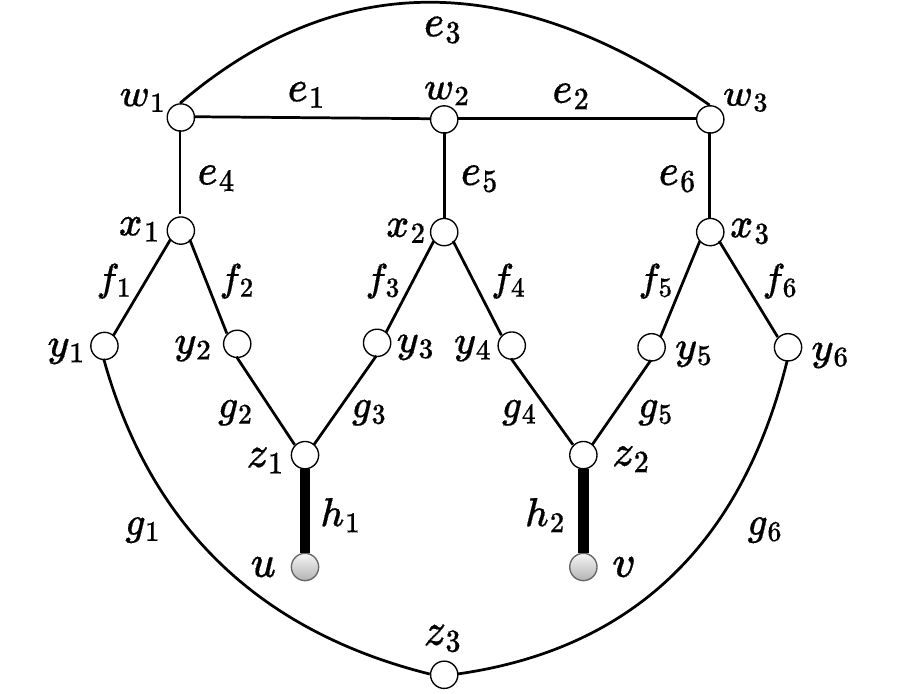}}
 	\caption{The graph $R[u,v]$.}
 	\label{fig:4}
 	\end{minipage}
 	\begin{minipage}[t]{7.9cm}
 		\centering
 		\resizebox{7.9cm}{6.1cm}{\includegraphics{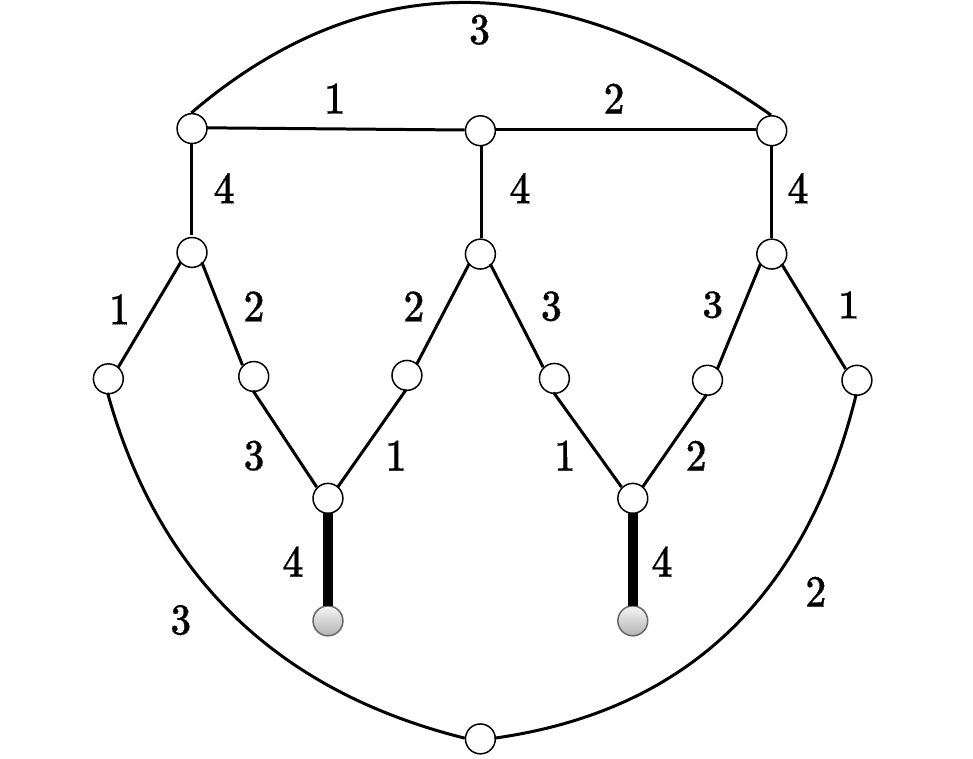}}
 	\caption{A semstrong 4-edge-coloring of $R[u,v]$.}
 	\label{fig:coloring}
 	\end{minipage}
 \end{figure}

Assume that  $G$ has a   $4$-edge-coloring 
$\phi$.
We color each subgraph $R[u,v]$ of $H$ corresponding to the edge $uv\in E(G)$ independently according to $\phi(uv)$ as follows.
Without loss of generality, we may assume that $\phi(uv)=4$. 
We first color the five edges $h_{1},h_{2},e_{4},e_{5},e_{6}$ with the same color 4. 
And then we color  the five edges $e_{1},g_{3},g_{4},f_{6},f_{1}$ with the same color 1,   the five  edges $e_{2},g_{5},g_{6},f_{2},f_{3}$ with the same color 2 and  the five   edges $e_{3},g_{1},g_{2},f_{4},f_{5}$ with the same color 3. Refer to \Cref{fig:coloring} for the coloring of $R[u,v]$.  This gives rise to a  $4$-edge-coloring of $H$ called  $\psi$.  It is easy to check  that    $\psi$ is semistrong.

Conversely, let  $\psi$ be a semistrong  $4$-edge-coloring of  $H$.
Let $uv$ be an edge in $G$ and   $R[u,v]$ be its corresponding subgraph of $H$. 
In the following, we are dedicated to proving $\psi(h_{1})=\psi(h_{2})$ and thus we can color $uv$ with the color $\psi(h_{1})$ to obtain a 4-edge-coloring of $G$.

Notice that $C_{3}=w_{1}w_{2}w_{3}$ is a 3-cycle in $R[u,v]$,
we must have	$|\psi(\cup_{i=1}^{6}\{e_{i}\})|=4$, $\psi(e_{4})=\psi(e_{5})=\psi(e_{6})$ and $\psi(e_{1}),\psi(e_{2}),\psi(e_{3})$ are different from each other.
Therefore, we may assume that $\psi(e_{1})=1$, $\psi(e_{2})=2$, $\psi(e_{3})=3$ and $\psi(e_{4})=\psi(e_{5})=\psi(e_{6})=4$.
Observe that the graph $R[u,v]$ has three  induced 7-cycles $C_{7}^{1}=z_{1}y_{3}x_{2}w_{2}w_{1}x_{1}y_{2}z_{1}$, $C_{7}^{2}=z_{2}y_{5}x_{3}w_{3}w_{2}x_{2}y_{4}z_{2}$ and 
$C_{7}^{3}=z_{3}y_{1}x_{1}w_{1}w_{3}x_{3}y_{6}z_{3}$. 
Recall that any 7-cycle has semistrong chromatic index 4.
Now  the color 4 occurs on each induced 7-cycle twice, we immediately observe the following.

\begin{observation}\label{obser:Quv-f-g}
	For each $i\in[1,6]$, $\psi(f_{i})\neq4$ and  $\psi(g_{i})\neq4$.
\end{observation}

In the following, we  prove   that  $\psi(h_{1})=\psi(h_{2})=4$. 
If not, by symmetry, we may assume that $\psi(h_{1})=3$. It  follows that $\psi(g_{2})\neq3$. 
Since  $\psi(h_{1})=\psi(e_{3})=3$ and both $h_{1}$ and $e_{3}$ are 2-neighbors of $f_{2}$, $\psi(f_{2})\neq3$. 
According to  Observation  \ref{obser:Quv-f-g}, $\psi(f_{2})\neq4$ and  $\psi(g_{2})\neq4$.
Therefore, we have $\{\psi(f_{2}),\psi(g_{2})\}=\{1,2\}$.

Now, if  $\psi(f_{2})=1$ and $\psi(g_{2})=2$, then $\psi(g_{3})\neq1$
since both $e_{1}$ and $g_{3}$ are 2-neighbors of $f_{2}$ and  $\psi(f_{2})=\psi(e_{1})=1$.
Because  $\psi(g_{2})=2$ and $\psi(h_{1})=3$, we must have $\psi(g_{3})=4$, which contradicts Observation \ref{obser:Quv-f-g}.

And if  $\psi(f_{2})=2$ and $\psi(g_{2})=1$, then by Observation  \ref{obser:Quv-f-g},
we have $\psi(g_{3})=2$ and thus $\psi(f_{3})\neq 2$.
Since $\psi(e_{1})=\psi(g_{2})=1$ and
 both $e_{1}$ and $g_{2}$ are 2-neighbors of $f_{3}$, $\psi(f_{3})\neq1$. 
 And  by Observation \ref{obser:Quv-f-g},  $\psi(f_{3})\neq 4$.
 Therefore, we have  $\psi(f_{3})=3$.
 It follows that  $\psi(f_{4})\in\{1,2\}$.
 Since $g_{3}$ has a 2-neighbor $f_{2}$ with the same color 2, $\psi(f_{4})\neq2$ and so $\psi(f_{4})=1$.
Then it follows that  $\psi(g_{4})\in\{2,3\}$.
We must have $\psi(g_{4})=2$ as otherwise  $f_{3}$ has no 1-vertex under $\psi$.
Moreover, due to $\psi(e_{1})=\psi(f_{4})=1$, $\psi(g_{4})=2$ and Observation \ref{obser:Quv-f-g}, we have $\psi(g_{5})=3$.
At this time, $f_{5}$ has two 2-neighbors $e_{2}$ and $g_{4}$ colored with 2, $\psi(f_{5})\neq2$.
Due to $\psi(g_{5})=3$  and Observation \ref{obser:Quv-f-g},  $\psi(f_{5})=1$.
It follows that $\psi(f_{6})\in\{2,3\}$.
Notice that  $e_{1}$ has a 2-neighbor $f_{4}$ with the same color 1, 
$\psi(f_{1})\neq1$ and thus $\psi(f_{1})=3$. 
This implies that $\psi(g_{1})\in\{1,2\}$.
We must have $\psi(g_{1})=1$ as otherwise $f_{2}$ has no 1-vertex under $\psi$.
Moreover, since   $\psi(f_{1})=3$, we must have $\psi(f_{6})=2$ as otherwise   $e_{3}$ has no 1-vertex under $\psi$.
It follows from $\psi(f_{6})=2$ and $\psi(g_{6})\neq4$ that 
$\psi(g_{6})=3$. At this time, $f_{1}$ has two $2$-neighbors $e_{3}$ and $g_{6}$ with the same color 3, and thus $f_{1}$  has no 1-vertex under $\psi$.
This gives a contradiction to the semistrong edge coloring $\psi$ of $H$.

Therefore, it holds that $\psi(h_{1})=\psi(h_{2})=4$.   For each edge $uv\in E(G)$, let $\phi(uv)=\psi(h_{1})=\psi(h_{2})$, where $h_{1}$ and $h_{2}$ are the edges of its corresponding subgraph $R[u,v]$ of $H$. We obtain a  $4$-edge-coloring $\phi$ of $G$.
 This completes  the proof of Case 2.

For each	 $k\ge3$ and each $k$-regular graph $G$, we  construct a corresponding  graph $H$ with maximum degree $k$ such that $G$ has a $k$-edge-coloring if and only if   $H$ has a semistrong $k$-edge-coloring.
Theorem \ref{Main-th} is proved.
\end{proof}

\section{The semistrong edge coloring of trees} \label{sec:tree}
This section designs a dynamic programming algorithm to solve the semistrong edge coloring problem for trees (\Cref{Main-th-tree}).  
Recall that any tree $T$ with maximum degree $\Delta$ is either semistrongly $\Delta$-edge-colorable or semistrongly $(\Delta+1)$-edge-colorable (see \Cref{sec:intro}), it suffices to  consider the following decision version of the semistrong edge coloring problem on trees.
    
\medskip	

\noindent{\textbf{Semistrong Edge Coloring Problem of Trees}}
\\
{\em Instance}: A tree $T$ with maximum degree $\Delta$.\\
{\em Question}: Does $T$ admit a semistrong $\Delta$-edge-coloring? 
    
\medskip	

\begin{figure}[H]
	\centering
	\resizebox{15.5cm}{5cm}{\includegraphics{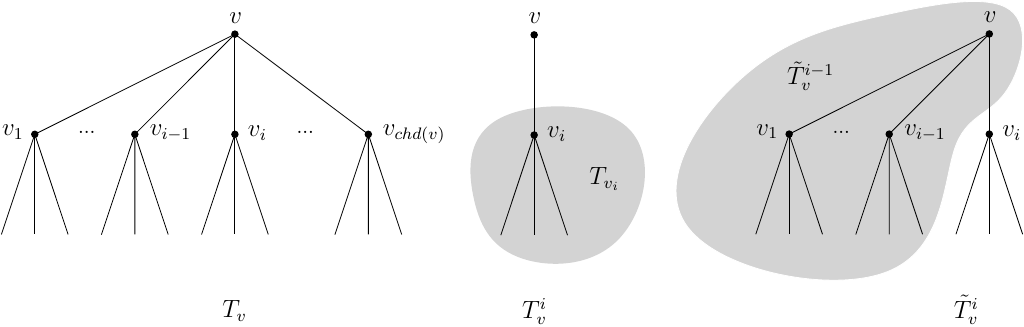}}
	\caption{The illustration of various subtrees of $T_r$.}
	\label{fig:Tv}
\end{figure}

Let $T_{r}$  be a tree rooted at a vertex $r$. If a vertex $v$ ($\neq r$) is of degree one in  $T_{r}$, then we call $v$ a {\em leaf} of $T_{r}$.
 For any $v\in V(T_{r})$, let $T_{v}$ denote the subtree of $T_{r}$ consisting of $v$ and all its descendants, and we denote by 
  $CHD(v)$ the set of children of $v$ in $T_{r}$  and $chd(v)$ the number of children of $v$ in $T_{r}$.
 For any $v\in V(T_{r})$,  let $v_{1},v_{2},\ldots,v_{chd(v)}$ denote the children of $v$. For $1\le i\le chd(v)$, let $T_{v}^{i}=T_{v_{i}}\cup\{vv_{i}\}$. Let $\tilde{T}_{v}^{1}=T_{v}^1$ and let $\tilde{T}_{v}^{i}=\tilde{T}_{v}^{i-1}\cup T_{v}^{i}$ for each $2\le i\le chd(v)$ if $chd(v)\ge2$. (See \Cref{fig:Tv} for an illustration.)
 It is clear that $T_{v}=\tilde{T}_{v}^{chd(v)}$.

 For any two vertices $x,y\in V(T_{v})$, let $d_{T_{v}}(x,y)$ denote  the  distance between  $x$ and $y$ in $T_{v}$. 
 For an edge $xy\in E(T_{v})$, define $d_{T_{v}}(xy,v)=\max\{d_{T_{v}}(x,v),d_{T_{v}}(y,v)\}$.
 If $d_{T_{v}}(xy,v)=i$, then we call $xy$ an {\em $i$th-generation edge descended from $v$} (we usually abbreviate it as {\em $i$th-generation edge}  in the absence of confusion).

Assume that   $T_{v}$ is semistrongly  $\Delta$-edge-colorable and let  $\phi$ be a 
 semistrong $\Delta$-edge-coloring  of $T_{v}$. We classify the $\Delta$ colors into the following five types with respect to $\phi$, according to whether a color appears on $1$st-generation edges or  $2$nd-generation edges descended from $v$ and, when they occur, the position of the  1-vertices of the corresponding edges.
\begin{definition}[Category of colors]
\label{def:colors}
Let  $\phi$ be a 
 semistrong $\Delta$-edge-coloring  of $T_{v}$, and let  $\alpha$ be a color under $\phi$. 
     \begin{itemize}
 \setlength\itemsep{2pt}
    \item   
    If $\alpha$  appears  on neither $1$st-generation edges nor  $2$nd-generation edges descended from $v$,  then  $\alpha$ is {\em a color of type $\mathcal{A}$ under $\phi$}.
     \item  If there is a  $1$st-generation edge $vw$ descended from $v$ such that $\phi(vw)=\alpha$ and $w$ is the $1$-vertex of $vw$ under $\phi$, then we say that $\alpha$ is {\em a color of type $\mathcal{P}$ under $\phi$}. (Note that in this case, none of the  $1$st- and $2$nd-generation edges descended from $w$ are  colored $\alpha$ under $\phi$.)
   \item  
  If there is a  $1$st-generation edge $vw$ descended from $v$ such that $\phi(vw)=\alpha$ and $w$ is not the 1-vertex of $vw$ under $\phi$ (hence $v$ must be the  1-vertex of $vw$ under $\phi$), then we say that $\alpha$ is  {\em a color of type $\mathcal{Q}$ under $\phi$}.
  \item  
 If $\alpha$  appears on $2$nd-generation edges but not $1$st-generation edges descended from $v$, and for every $2$nd-generation edge $xy$ colored $\alpha$ with $d_{T_{v}}(x,v)<d_{T_{v}}(y,v)$, $y$  is the 1-vertex of $xy$ under $\phi$, then we call  $\alpha$  {\em the color of type $\mathcal{S}$ under $\phi$}.
   \item  
  If $\alpha$  appears on $2$nd-generation edges but not $1$st-generation edges descended from $v$, and there is a $2$nd-generation edge $xy$ colored $\alpha$ with $d_{T_{v}}(x,v)<d_{T_{v}}(y,v)$ such that $y$  is not the 1-vertex of $xy$ under $\phi$ (hence  $x$ must be the  1-vertex of $xy$ under $\phi$), then we call  $\alpha$  {\em the color of type $\mathcal{T}$ under $\phi$}.
 \end{itemize}

 We use  $\mathcal{A}_{\phi}$ (resp. $\mathcal{P}_{\phi}$, $\mathcal{Q}_{\phi}$, $\mathcal{S}_{\phi}$, and  $\mathcal{T}_{\phi}$) to denote the set of colors of type 
$\mathcal{A}$ (resp. $\mathcal{P}$, $\mathcal{Q}$, $\mathcal{S}$ and  $\mathcal{T}$) under $\phi$.
\end{definition}

\begin{figure}[H]
	\centering
	\resizebox{15.8cm}{7.8cm}{\includegraphics{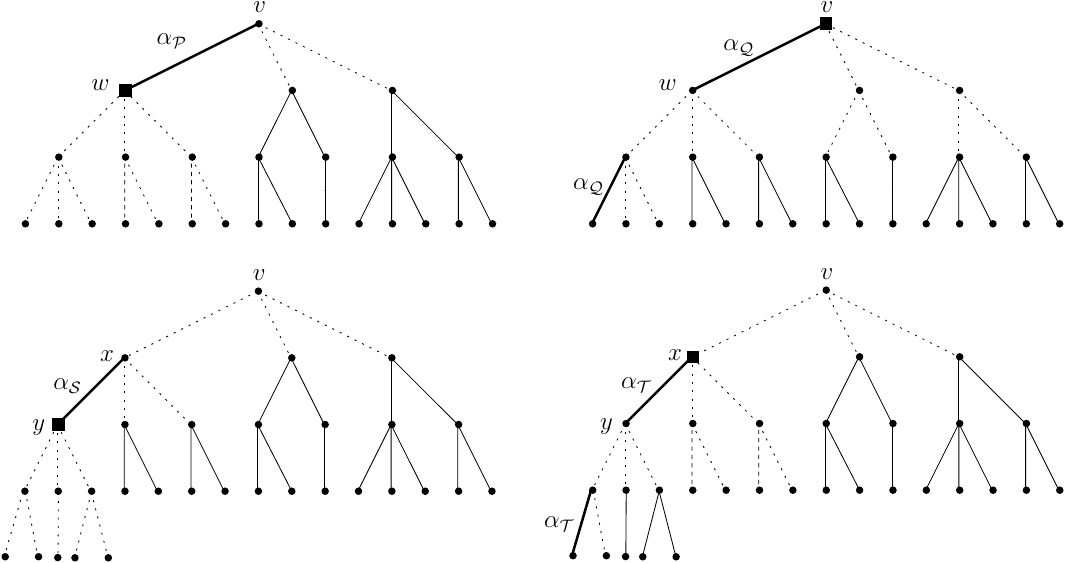}}
    \caption{
   Colors of Type $\mathcal{P}$, $\mathcal{Q}$, $\mathcal{S}$ and $\mathcal{T}$  under a semistrong $\Delta$-coloring of $T_v$.
    Black square dots indicate vertices that must be the 1-vertices of the corresponding colored edges, and dashed edges indicate edges that cannot be colored with that color.
   }
	\label{fig:color}
\end{figure}

   See \Cref{fig:color} for an illustration of the classification of colors.  
    Note that if a color $\alpha$ is of Type $\mathcal{T}$, then it is possible that there is one $2$nd-generation edge $x'y'$ colored $\alpha$ with $d_{T_{v}}(x',v)<d_{T_{v}}(y',v)$ such that $y'$  is  the 1-vertex of $xy$ under $\phi$.

Given four nonnegative integers  $p,q,s$ and $t$ with $p+q+s+t\le \Delta$. 
A {\em $(p,q;s,t)$-coloring} of $T$  is a semistrong $\Delta$-edge-coloring with exactly $p,q,s$ and $t$ colors of type $\mathcal{P}$  $\mathcal{Q}$, $\mathcal{S}$ and $\mathcal{T}$, respectively (hence $|\mathcal{A}|=\Delta-p-q-s-t$). 
We now introduce the Boolean variable $b(T;p,q;s,t)$ to indicate whether the tree $T$ has  a $(p,q;s,t)$-coloring:
\[b(T; p,q;s,t) := \begin{cases} 1, & \text{if } T \text{ has a } (p,q;s,t)\text{-coloring},\\ 0, & \text{otherwise}. \end{cases} \]
We denote by $\mathcal{B}(T)$ the set of all tuples $(p,q;s,t)$ such that $T$ has a $(p,q;s,t)$-coloring, that is,
$$\mathcal{B}(T):=\{(p,q;s,t):b(T;p,q;s,t)=1\}.$$
For any two nonnegative integers $p_{0}$  and $q_{0}$, we define $$\mathcal{B}(T)|_{(p_{0},q_{0})}:=\{(p_{0},q_{0};s,t):(p_{0},q_{0};s,t)\in \mathcal{B}(T)\},$$ 
that is, the subset of $\mathcal{B}(T)$  consisting of all tuples whose first two coordinates are fixed to be $(p_{0},q_{0})$.

The next three useful facts immediate from the above definitions.

\begin{fact}\label{fact:colorable}
 Let $T$ be a tree with maximum degree $\Delta$. Then, $\chi_{ss}'(T)=\Delta$ if and only if $\mathcal{B}(T)\neq\emptyset$.
\end{fact}

\begin{fact}\label{fact:leaf}
  Let $v$ be a leaf of  a tree $T$. Then,  $\mathcal{B}(T)=\{(0,0;0,0)\}$.
\end{fact}

\begin{fact}\label{fact:empty}
  Let $T'$ be a subtree of the tree $T$ (with maximum degree $\Delta$).
  If $\mathcal{B}(T')=\emptyset$, then $\mathcal{B}(T_{r})=\emptyset$
(and thus $\chi_{ss}'(T_{r})=\Delta+1$).
\end{fact}

According to \Cref{fact:colorable}, determining whether $T_r$ is  semistrongly $\Delta$-edge-colorable reduces to checking whether $\mathcal{B}(T_{r})$ is empty.  The computation of 
$\mathcal{B}(T_{r})$
 is based on a dynamic programming approach, as described below.
Firstly, we assign initial values  $\mathcal{B}(T_{v})$ to all leaves  (using \Cref{fact:leaf}) in $T_r$.
Secondly, for a non-leaf vertex $v$ we compute  $\mathcal{B}(T_{v})$ based on $\{\mathcal{B}(T_{v_i}):v_i \in CHD(v)\}$ in two stages:
\begin{itemize}
 \setlength\itemsep{2pt}
    \item (Vertical expansion) For each child $v_i$ of $v$, compute $\mathcal{B}(T_{v}^{i})$ based on $\mathcal{B}(T_{v_{i}})$ in $O(\Delta^{4})$ time, as presented in  Lemma~\ref{lemma:vertical}. See Subsection~\ref{subsec:vertical} for details.
     \item (Horizontal merging) For $i$ from $1$ to $chd(v)-1$, compute  $\mathcal{B}(\tilde{T}_{v}^{i+1})$  based on    $\mathcal{B}(\tilde{T}_{v}^{i})$  and  $\mathcal{B}(T_{v}^{i+1})$ in $O(\Delta^{6})$ time, as presented in   Lemma~\ref{lemma:horizontal}.  See Subsection~\ref{subsec:horizon} for details.
\end{itemize}
If we meet a vertex $v$ satisfying $\mathcal{B}(T_{v})=\emptyset$, then by \Cref{fact:empty} we have   $\chi_{ss}'(T_{r})=\Delta+1$; otherwise,$\chi_{ss}'(T_{r})=\Delta$.
The  dynamic programming algorithm is presented below.

\makeatletter
\newenvironment{breakablealgorithm}
{
		\begin{center}
			\refstepcounter{algorithm}
			\hrule height.8pt depth0pt \kern2pt
			\renewcommand{\caption}[2][\relax]{
				{\raggedright\textbf{\ALG@name~\thealgorithm} ##2\par}%
				\ifx\relax##1\relax 
				\addcontentsline{loa}{algorithm}{\protect\numberline{\thealgorithm}##2}%
				\else 
				\addcontentsline{loa}{algorithm}{\protect\numberline{\thealgorithm}##1}%
				\fi
				\kern2pt\hrule\kern2pt
			}
		}{
		\kern2pt\hrule\relax
	\end{center}
}
\makeatother

\begin{breakablealgorithm} 
	\caption{ Semistrong edge coloring of trees}  
	\label{alg:main}  
	\begin{algorithmic}[1]  
		\Require   a tree $T_{r}$ with maximum degree $\Delta$;	
		\Ensure $\chi_{ss}'(T_{r})$;
		
		\While{there is an uncolored leaf $v$}
		\State
		set 	 		$\mathcal{B}(T_{v}):=\{(0,0;0,0)\}$
		\State color $v$  red	
		\EndWhile	
		\\
		
		\While{there is an uncolored vertex $v$ whose children are all colored red}	
	
		\For{$i:1\to chd(v)$}
		\State	
		call Algorithm \ref{alg:vertical}   on $\mathcal{B}(T_{v_{i}})$ to compute $\mathcal{B}(T_{v}^{i})$
		
		\If{$\mathcal{B}(T_{v}^{i})=\emptyset$}
			\State 	stop and output $\chi_{ss}'(T_{r})=\Delta+1$
			\EndIf
		
		\EndFor

			\For{$i:1\to chd(v)-1$}
		\State	
		call Algorithm \ref{alg:horizontal} on $\mathcal{B}(\tilde{T}_{v}^{i})$ and $\mathcal{B}(T_{v}^{i+1})$  to compute $\mathcal{B}(\tilde{T}_{v}^{i+1})$
		
		\If{$\mathcal{B}(\tilde{T}_{v}^{i+1})=\emptyset$}
		\State	stop and output  $\chi_{ss}'(T_{r})=\Delta+1$
		\EndIf
		
		\EndFor
		
		\State Color $v$  red		
		
	\EndWhile
    \State	output  $\chi_{ss}'(T_{r})=\Delta$
	

	\end{algorithmic}  
\end{breakablealgorithm} 

Before proceeding to the next subsection, we now analyze the time complexity of the above algorithm.
Let $v$ be a vertex of $T_r$. Suppose that $\{\mathcal{B}(T_{v_i}):v_i\in CHD(v)\}$ is known. One can first compute  $\{\mathcal{B}(T_{v}^{i}) :v_i\in CHD(v)\}$  in $O(\Delta^{4}d_{T_r}(v))$ time (by Lemma~\ref{lemma:vertical}) and then further 
obtain   $\mathcal{B}(T_{v})$  in $O(\Delta^{6}d_{T_r}(v))$ time (by Lemma~\ref{lemma:horizontal}).
It follows that 
\Cref{alg:main} can compute  $\mathcal{B}(T_{r})$ or output $\chi_{ss}'(T_{r})=\Delta+1$
in $O(\Delta^{6}\sum_{v\in V(T_r)}d_{T_r}(v))=O(\Delta^{6}n)$ time. 
Hence, one
can determine the semistrong chromatic index of $T_r$ in  $O(\Delta^{6}n)$ time and so \Cref{Main-th-tree}   follows.

\medskip
In the subsequent two subsections, we prove Lemmas~\ref{lemma:vertical} and~\ref{lemma:horizontal}, respectively.


\subsection{Vertical expansion}
\label{subsec:vertical}
In this subsection, we focus on computing $\mathcal{B}(T_{v}^{i})$ based on  $\mathcal{B}(T_{v_{i}})$ ($\neq \emptyset$).  
Since $T_{v}^{i}=T_{v_{i}}\cup\{vv_{i}\}$ (see \Cref{fig:Tv}), 
we have $\mathcal{B}(T_{v}^{i})=\mathcal{B}(T_{v}^{i})|_{(1,0)}\cup\mathcal{B}(T_{v}^{i})|_{(0,1)}$.  Note that
\begin{itemize}
\setlength\itemsep{2pt}
    \item  $\mathcal{B}(T_{v}^{i})|_{(1,0)}$ indicates that  whether there exists a semistrong $\Delta$-edge-coloring of $T_{v}^{i}$ such that 
$v_i$ is a 1-vertex of $vv_i$. 
\item $\mathcal{B}(T_{v}^{i})|_{(0,1)}$ indicates that
whether there exists a semistrong $\Delta$-edge-coloring of $T_{v}^{i}$ such that 
$v_i$ is not a 1-vertex of $vv_i$ and $v$ is a 1-vertex of $vv_i$.
\end{itemize}
Therefore, the problem of computing $\mathcal{B}(T_{v}^{i})$ can be simplified by separately computing $\mathcal{B}(T_{v}^{i})|_{(1,0)}$  and $\mathcal{B}(T_{v}^{i})|_{(0,1)}$.
The idea to compute  $\mathcal{B}(T_{v}^{i})|_{(1,0)}$ (resp. $\mathcal{B}(T_{v}^{i})|_{(1,0)}$) is to determine whether there exists a  $(p,q;s,t)$-coloring of $T_{v_i}$ that  can be extended into a $(1,0;p,q)$-coloring (resp. $(0,1;p,q)$-coloring) of $T_{v}^i$. 
Towards to this, we first prove the following claim.

 \begin{claim}
Let $\phi$ be    a  $(p,q;s,t)$-coloring of $T_{v_i}$. Then,
 \begin{itemize}
     \item $\phi$ can be extended to  a $(1,0;p,q)$-coloring of $T_{v}^i$ if and only if $p+q+s+t\le \Delta-1$;
     \item $\phi$  can be extended to  a $(0,1;p,q)$-coloring of $T_{v}^i$ if and only if $p+q+s+t\le \Delta$ and $s\ge1$.
 \end{itemize}
 \end{claim}
 \begin{proof}
First assume  that     $\phi$ can be extended to  a $(1,0;p,q)$-coloring  $\psi$ of $T_{v}^i$. 
Then  none of the $1$st- and $2$nd-generation edges descended from $v_i$ in  $T_{v_i}$  is colored $\psi(vv_i)$ under $\phi$. 
Hence,  $\psi(vv_i)$ is of type $\mathcal{A}$ (see Definition~\ref{def:colors}) under $\phi$, and so  $p+q+s+t\le \Delta-1$.
Conversely, assume $p+q+s+t\le \Delta-1$. Then there exists a color $\alpha$ that does not appear on any of the $1$st- and $2$nd-generation edges descended from $v_i$ in  $T_{v_i}$ under $\phi$. It is  easy to extend $\phi$  to  a $(1,0;p,q)$-coloring  $\psi$ of $T_{v}^i$ by coloring $vv_i$ with $\alpha$.

Now, assume that $\phi$  can be extended to  a $(0,1;p,q)$-coloring $\psi$ of $T_{v}^i$. Then, none of the $1$st-generation edges descended from $v_i$ in  $T_{v_i}$  are colored $\phi(vv_i)$,  and  some $2$nd-generation edges descended from $v_i$ in  $T_{v_i}$  are colored $\phi(vv_i)$ under $\phi$.  Thus,
     the color $\phi(vv_i)$ is of  type $\mathcal{S}$ (see Definition~\ref{def:colors}) under $\phi$, and so $p+q+s+t\le \Delta$ and $s\ge1$.
  Conversely, assume $p+q+s+t\le \Delta$ and $s\ge1$. 
     Then there exists a color $\alpha$ of type $\mathcal{S}$ under $\phi$. 
     By coloring $vv_i$ with $\alpha$, 
     we extend $\phi$  to  a $(1,0;p,q)$-coloring  $\psi$ of $T_{v}^i$.
 \end{proof}
Based on the above claim, we
now present our algorithm.


\begin{breakablealgorithm} 
	\caption{ Compute $\mathcal{B}(T_{v}^{i})$ from  $\mathcal{B}(T_{v_{i}})$}  
	\label{alg:vertical}  
	\begin{algorithmic}[1]  
		\Require  $\mathcal{B}(T_{v_{i}})$ ($\neq \emptyset$);	
		\Ensure $\mathcal{B}(T_{v}^{i})$;
		\State
	set 		$\mathcal{B}(T_{v}^{i}):=\emptyset$
	
	\For{$p:0\to chd(v_{i})$}
		\State
set 	 $q:=chd(v_{i})-p$

\If{there are $s,t\in \mathbb{N}$ such that $(p,q;s,t)\in \mathcal{B}(T_{v_{i}})$ and $p+q+s+t\le \Delta-1$}

\State
set 		$\mathcal{B}(T_{v}^{i}):=\mathcal{B}(T_{v}^{i})\cup\{(1,0;p,q)\}$

\EndIf

\If{there are $s,t\in \mathbb{N}$ such that $(p,q;s,t)\in \mathcal{B}(T_{v_{i}})$, $p+q+s+t\le \Delta$ and $s\ge1$}

\State
set 		$\mathcal{B}(T_{v}^{i}):=\mathcal{B}(T_{v}^{i})\cup\{(0,1;p,q)\}$

\EndIf
	\EndFor

	\State	
\Return	 $\mathcal{B}(T_{v}^{i})$

	\end{algorithmic}  
\end{breakablealgorithm}

The following lemma is easy to see.
\begin{lemma}\label{lemma:vertical}
Given a vertex $v\in V(T_{r})$ and one of its children $v_{i}$. Then, Algorithm \ref{alg:vertical} can compute  $\mathcal{B}(T_{v}^{i})$   based on $\mathcal{B}(T_{v_{i}})$
in $O(\Delta^{3})$ time.
\end{lemma}

\subsection{Horizontal merging}
	\label{subsec:horizon}
In this subsection, we are dedicated to computing  $\mathcal{B}(\tilde{T}_{v}^{i+1})$  from   $\mathcal{B}(\tilde{T}_{v}^{i})$  and  $\mathcal{B}(T_{v}^{i+1})$.
The main result of this subsection is presented in
Lemma~\ref{lemma:horizontal}.

The crucial step in calculating $\mathcal{B}(\tilde{T}_{v}^{i+1})$ lies in answering the following question: Given four nonnegative integers $p,q,s,t$, is there a  $(p,q;s,t)$-coloring $\phi$ of $\tilde{T}_{v}^{i+1}$?
 To address this question, we establish the following lemma, which characterizes precisely when a 
$(p,q;s,t)$-coloring exists.

\begin{lemma}\label{lemma:iff}
Let $p,q,s,t$  be  four nonnegative integers  with $p+q+s+t\le\Delta$. Then, there exists a  $(p,q;s,t)$-coloring of $\tilde{T}_{v}^{i+1}$ 
if and only if  there exist four nonnegative integers $s_L,t_L,s_R,t_R$
such that
at least one of the following two cases holds.
  \begin{itemize}
    \setlength\itemsep{2pt}
        \item $(p,q-1;s_{L},t_{L})\in\mathcal{B}(\tilde{T}_{v}^{i})$,  $(0,1;s_{R},t_{R})\in\mathcal{B}(T_{v}^{i+1})$, and
        \[
\quad(1.1)^*
\left\{\!
\begin{array}{r@{\hspace{0.5em}}ll}
(1).&t_R\le t \le t_L+t_R, \\
(2).& \max\{0,s_L-s\}\le \min\{s_L,t-t_L\},\\
(3).& 0\le s_L+s_R-s\le s_L+p+t-t_R.
\end{array}
\right.
\]
  \item $(p-1,q;s_{L},t_{L})\in\mathcal{B}(\tilde{T}_{v}^{i})$,  $(1,0;s_{R},t_{R})\in\mathcal{B}(T_{v}^{i+1})$, and either
  \[
\quad(1.2)^*
\left\{\!
\begin{array}{r@{\hspace{0.5em}}ll}
(1).&t_R\le t \le t_L+t_R, \\
(2).& \max\{0,s_L-s\}\le \min\{s_L,t-t_L\},\\
(3).& 0\le s_L+s_R-s\le s_L+(p-1)+t-t_R.
\end{array}
\right.
\]
or \[
\quad(1.3)^*
\left\{\!
\begin{array}{r@{\hspace{0.5em}}ll}
(1).&t_R\le t \le t_L+t_R, \\
(2).& \max\{0,(s_L-1)-s\}\le \min\{s_L-1,t-t_L\},\\
(3).& 0\le (s_L-1)+s_R-s\le (s_L-1)+(p-1)+t-t_R.
\end{array}
\right.
\]
    \end{itemize}
   $($Note that  by regarding $p$
 in  $(1.1)^*$ as $p-1$, we obtain  $(1,2)^*$; 
 by replacing  $p$ as $p-1$ and $s_L$ as $s_L-1$ in $(1.1)^*$, we have $(1.3)^*$.$)$
\end{lemma}

Before proving this lemma,
we first present the algorithm for computing $\mathcal{B}(\tilde{T}_{v}^{i+1})$, whose design is grounded in the lemma.

\begin{breakablealgorithm} 
	\caption{Compute  $\mathcal{B}(\tilde{T}_{v}^{i+1})$  from   $\mathcal{B}(\tilde{T}_{v}^{i})$  and  $\mathcal{B}(T_{v}^{i+1})$}  
	\label{alg:horizontal}  
	\begin{algorithmic}[1]  
		\Require  $\mathcal{B}(\tilde{T}_{v}^{i})$ ($\neq\emptyset$) and  $\mathcal{B}(T_{v}^{i+1})$ ($\neq\emptyset$);	
		\Ensure $\mathcal{B}(\tilde{T}_{v}^{i+1})$;
		\State
		set 		$\mathcal{B}(\tilde{T}_{v}^{i+1}):=\emptyset$
		
		\For{$p:0\to i+1$}
			\State
		set  $q:=(i+1)-p$
		
		\For{$s:0\to \Delta-p-q$}
	\For{$t:0\to \Delta-p-q-s$}

					\If{there exists a  $(p,q;s,t)$-coloring of $\tilde{T}_{v}^{i+1}$ (verified using Lemma~\ref{lemma:iff})}	
				\State set  $\mathcal{B}(\tilde{T}_{v}^{i+1}):=\mathcal{B}(\tilde{T}_{v}^{i+1})\cup\{(p,q;s,t)\}$
		
			\EndIf

	\EndFor
	\EndFor		
		\EndFor		
		
		\State	
	\Return $\mathcal{B}(\tilde{T}_{v}^{i+1})$
				
	\end{algorithmic}  
\end{breakablealgorithm} 

According to Lemma~\ref{lemma:iff} and  \Cref{alg:horizontal}, 
one can directly  deduce  the following lemma. 

 \begin{lemma}\label{lemma:horizontal}
	Given a vertex $v\in V(T_{r})$. For each $1\le i\le chd(v)-1$, \Cref{alg:horizontal} can compute  $\mathcal{B}(\tilde{T}_{v}^{i+1})$ based on 
  $\mathcal{B}(\tilde{T}_{v}^{i})$ and  $\mathcal{B}(T_{v}^{i+1})$  in $O(\Delta^{6})$ time.
\end{lemma}
\begin{proof}
Given $p,q,s,t\ge0$. 
Note first that if $(0,1;s_{R},t_{R})\in\mathcal{B}(T_{v}^{i+1})$ or  $(1,0;s_{R},t_{R})\in\mathcal{B}(T_{v}^{i+1})$ then we must have 
$s_{R}+t_{R}=chd(v_{i+1})$.
It implies that one can determine whether there  exist four nonnegative integers $s_L,t_L,s_R,t_R$ satisfying Lemma~\ref{lemma:iff} (and thus whether there is  a  $(p,q;s,t)$-coloring  of $\tilde{T}_{v}^{i+1}$)  in $O(\Delta^3)$ time.
Moreover, if  $\tilde{T}_{v}^{i+1}$ has a  $(p,q;s,t)$-coloring  of $\tilde{T}_{v}^{i+1}$, then  $p+q=i+1$.
 Hence,  one can easily  
obtain  $\mathcal{B}(\tilde{T}_{v}^{i+1})$  in $O(\Delta^{6})$ time by  \Cref{alg:horizontal}.
\end{proof}

 

\medskip
In the remainder of this subsection, we concentrate on proving Lemma~\ref{lemma:iff}.
\medskip

\noindent\textbf{\textit{Proof of Lemma~\ref{lemma:iff}.}}
Let $\phi$ be  a $(p,q;s,t)$-coloring  of $\tilde{T}_{v}^{i+1}$. 
We denote by $\phi_{L}$ the coloring  $\phi$ restricted to the subtree  $\tilde{T}_{v}^{i}$ and  $\phi_{R}$ the coloring  $\phi$ restricted to the subtree  $T_{v}^{i+1}$. 
Suppose that  $\phi_{L}$  is a   $(p_L,q_L;s_L,t_L)$-coloring  of $\tilde{T}_{v}^{i}$ and $\phi_{R}$ 
is a   $(p_R,q_R;s_R,t_R)$-coloring of   $T_{v}^{i+1}$.
 We have the following simple but important observation.

\begin{observation}
\label{obser:3cases}
	$\phi_{L}$ and $\phi_{R}$ satisfy exactly
	one of the following three situations:
    \begin{enumerate}
    \setlength\itemsep{2pt}
        \item $\phi_{L}$ is  a $(p,q-1;s_{L},t_{L})$-coloring of  $\tilde{T}_{v}^{i}$  and $\phi_{R}$ is  a  $(0,1;s_{R},t_{R})$-coloring of  $T_{v}^{i+1}$. $($Note that $\phi_{R}(vv_{i+1})\in\mathcal{A}_{\phi_{L}}$ at this time$)$.
        \item $\phi_{L}$ is  a $(p-1,q;s_{L},t_{L})$-coloring of  $\tilde{T}_{v}^{i}$  and $\phi_{R}$ is  a  $(1,0;s_{R},t_{R})$-coloring of  $T_{v}^{i+1}$ with $\phi_{R}(vv_{i+1})\in\mathcal{A}_{\phi_{L}}$.
        \item $\phi_{L}$ is  a $(p-1,q;s_{L},t_{L})$-coloring of  $\tilde{T}_{v}^{i}$  and $\phi_{R}$ is  a  $(1,0;s_{R},t_{R})$-coloring of  $T_{v}^{i+1}$ with $\phi_{R}(vv_{i+1})\in\mathcal{S}_{\phi_{L}}$.
     \end{enumerate}
\end{observation}

According to  $\phi_{L}$,
  the $\Delta$ colors  are divided into the five subsets   $\mathcal{P}_{\phi_{L}}$,  $\mathcal{Q}_{\phi_{L}}$, $\mathcal{S}_{\phi_{L}}$, $\mathcal{T}_{\phi_{L}}$ and  $\mathcal{A}_{\phi_{L}}$, where $|\mathcal{P}_{\phi_{L}}|=p_L$,  $|\mathcal{Q}_{\phi_{L}}|=q_L$, $|\mathcal{S}_{\phi_{L}}|=s_{L}$, $|\mathcal{T}_{\phi_{L}}|=t_{L}$ and  $|\mathcal{A}_{\phi_{L}}|=a_{L}=\Delta-p_L-q_L-s_{L}-t_{L}$.  
  (Note that we always have $a_{L}=\Delta+1-p-q-s_{L}-t_{L}$.)
  This partition is unique up to the permutation of colors.

 Let $\alpha=\phi_{R}(vv_{i+1})$. For $\mathcal{I}\in\{\mathcal{P},\mathcal{Q},\mathcal{S},\mathcal{T},\mathcal{A}\}$, let $\mathcal{X}_{\mathcal{I}}=\mathcal{S}_{\phi_{R}}\cap\mathcal{I}_{\phi_{L}}$ and $\mathcal{Y}_{\mathcal{I}}=\mathcal{T}_{\phi_{R}}\cap\mathcal{I}_{\phi_{L}}$, and let $x_{I}=|\mathcal{X}_{\mathcal{I}}|$ and $y_{I}=|\mathcal{Y}_{\mathcal{I}}|$.
According to the definitions of the five types of colors, we have $\mathcal{X}_{\mathcal{Q}}=\mathcal{Y}_{\mathcal{P}}=\mathcal{Y}_{\mathcal{Q}}=\emptyset$, $\mathcal{S}_{\phi_{R}}=\mathcal{X}_{\mathcal{P}}\cup \mathcal{X}_{\mathcal{S}}\cup \mathcal{X}_{\mathcal{T}}\cup \mathcal{X}_{\mathcal{A}}$ and 
 $\mathcal{T}_{\phi_{R}}=\mathcal{Y}_{\mathcal{S}}\cup \mathcal{Y}_{\mathcal{T}}\cup \mathcal{Y}_{\mathcal{A}}$.
 
\begin{figure}[H]  
	\centering
	\resizebox{13.5cm}{3cm}{\includegraphics{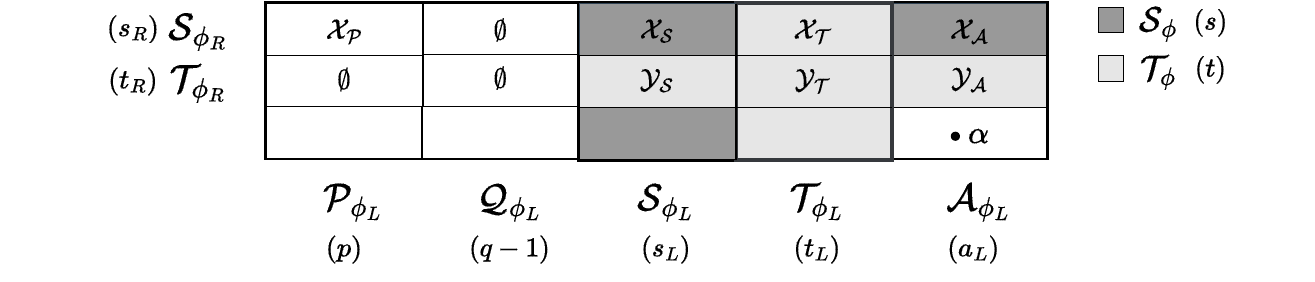}}
	\caption{The illustration of the first case in \Cref{obser:3cases}. (Note that $\alpha$ is not necessarily the only color present in the corresponding grid.)}
	\label{fig:Q1}
\end{figure}

Moreover, if  $\phi_{L}$ and $\phi_{R}$ satisfy the first case 
in \Cref{obser:3cases}, then   the following four properties hold. (See \Cref{fig:Q1} for an illustration.)

\begin{itemize}
\setlength\itemsep{2pt}
    \item $\mathcal{P}_{\phi}=\mathcal{P}_{\phi_{L}}$ (and  $|\mathcal{P}_{\phi}|=p$).
    \item $\mathcal{Q}_{\phi}=\mathcal{Q}_{\phi_{L}}\cup\{\alpha\}$ (and  $|\mathcal{Q}_{\phi}|=q$).
    \item  $\mathcal{S}_{\phi}=(\mathcal{S}_{\phi_{L}}\setminus \mathcal{Y}_{\mathcal{S}})\cup \mathcal{X}_{\mathcal{A}}$ and $|\mathcal{S}_{\phi}|=s$.
    \item $\mathcal{T}_{\phi}=\mathcal{T}_{\phi_{L}}\cup \mathcal{Y}_{\mathcal{S}}\cup \mathcal{Y}_{\mathcal{A}}$ and $|\mathcal{T}_{\phi}|=t$.
\end{itemize}

\begin{figure}[H]
	\centering
	\resizebox{13.5cm}{3cm}{\includegraphics{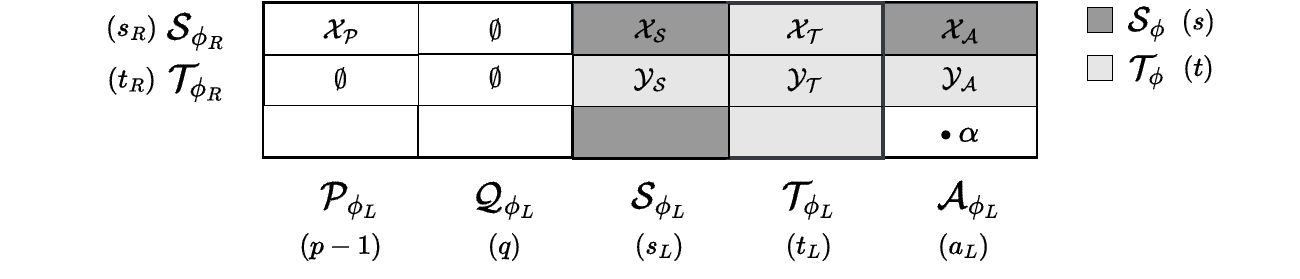}}
	\caption{The illustration of the second case in \Cref{obser:3cases}. (Note that $\alpha$ is not necessarily the only color present in the corresponding grid.)}
	\label{fig:Q2}
\end{figure}

\begin{figure}[htbp]  
	\centering
	\resizebox{13.5cm}{3cm}{\includegraphics{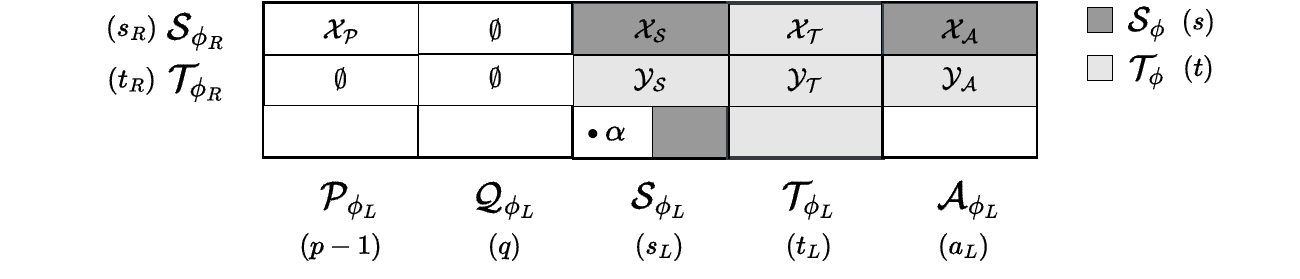}}
	\caption{The illustration of the third case in \Cref{obser:3cases}. (Note that $\alpha$ is the only color in the corresponding grid.)}
	\label{fig:Q3}
\end{figure}

Similarly, if $\phi_{L}$ and $\phi_{R}$ satisfy the second case (or the third case) 
in \Cref{obser:3cases}, analogous properties can be derived (details are omitted here), as illustrated in \Cref{fig:Q2} (or \Cref{fig:Q3}).

As illustrated in \Cref{fig:Q1,fig:Q2,fig:Q3}, for each of the three cases  in \Cref{obser:3cases},  the 
seven nonnegative integers $x_{P},x_{S},x_{T},x_{A},y_{S},y_{T},y_{A}$ satisfy the corresponding system of linear equalities and inequalities below, with each system derived directly from the respective figure and listed as $(1.1)$, $(1.2)$ and $(1.3)$.
\vspace{-1.0em}

\noindent
\begin{minipage}[t]{0.32\linewidth}
\small
\[
\quad(1.1)
\left\{\!\!
\begin{array}{r@{\hspace{0.35em}}l}
1.& s_{R} = x_{P} + x_{S} + x_{T} + x_{A}, \\
2.& t_{R} = y_{S} + y_{T} + y_{A}, \\
3.& s = s_{L} - y_{S} + x_{A}, \\
4.& t = t_{R} + t_{L} - y_{T},\\
5.& x_{P} \le p, \\
6.& x_{S} + y_{S} \le s_{L}, \\
7.& x_{T} + y_{T} \le t_{L}, \\
8.& 1 + x_{A} + y_{A} \le a_{L}.
\end{array}
\right.
\]
\end{minipage}%
\begin{minipage}[t]{0.32\linewidth}
\small
\[
\quad(1.2)
\left\{\!\!
\begin{array}{r@{\hspace{0.35em}}l}
1.& s_{R}=x_{P}+x_{S}+x_{T}+x_{A},\\
2.& t_{R}=y_{S}+y_{T}+y_{A}, \\
3. &  s=s_{L}-y_{S}+x_{A}, \\
4.& t=t_{R}+t_{L}-y_{T}, \\
 5.&   x_{P}\le p-1, \\
 6. &  x_{S}+y_{S}\le s_{L}, \\
7. & x_{T}+y_{T}\le t_{L}, \\
8. & 1+x_{A}+y_{A}\le a_{L}.
\end{array}
\right.
\]
\end{minipage}%
\begin{minipage}[t]{0.32\linewidth}
\small
\[\quad(1.3)
\left\{\!\!
\begin{array}{r@{\hspace{0.35em}}l}
1.& s_{R}=x_{P}+x_{S}+x_{T}+x_{A}, \\
2. & t_{R}=y_{S}+y_{T}+y_{A}, \\
3. & s=s_{L}-1-y_{S}+x_{A}, \\
4. & t=t_{R}+t_{L}-y_{T}, \\
5. & x_{P}\le p-1, \\
6.& 1+x_{S}+y_{S}\le s_{L}, \\
7. &  x_{T}+y_{T}\le t_{L}, \\
8. & x_{A}+y_{A}\le a_{L}.
\end{array}
\right.
\]
\end{minipage}

\bigskip
The next claim reduces the problem of determining whether $\tilde{T}_{v}^{i+1}$ has a $(p,q;s,t)$-coloring  to
 finding an integer solution to Systems $(1.1)$,  $(1.2)$ or $(1.3)$.
\begin{claim}
\label{claim:iff1}
 Given  eight nonnegative integers $p,q,s,t,s_{L},t_{L},s_{R},t_{R}$. Then, 
 \begin{itemize}
     \item  there exists a $(p,q;s,t)$-coloring  $\phi$  of $\tilde{T}_{v}^{i+1}$ such that
 $\phi_{L}$ is  a $(p,q-1;s_{L},t_{L})$-coloring of  $\tilde{T}_{v}^{i}$  and $\phi_{R}$ is  a  $(0,1;s_{R},t_{R})$-coloring of  $T_{v}^{i+1}$,
if and only if  $(p,q-1;s_{L},t_{L})\in\mathcal{B}(\tilde{T}_{v}^{i})$,  $(0,1;s_{R},t_{R})\in\mathcal{B}(T_{v}^{i+1})$, and
 there are seven nonnegative integers $x_{P},x_{S},x_{T},x_{A},y_{S},y_{T},y_{A}$  satisfying  System $(1.1)$.
 \item  there exists a $(p,q;s,t)$-coloring  $\phi$  of $\tilde{T}_{v}^{i+1}$ such that
 $\phi_{L}$ is  a $(p-1,q;s_{L},t_{L})$-coloring of  $\tilde{T}_{v}^{i}$  and $\phi_{R}$ is  a  $(1,0;s_{R},t_{R})$-coloring of  $T_{v}^{i+1}$,
if and only if  $(p-1,q;s_{L},t_{L})\in\mathcal{B}(\tilde{T}_{v}^{i})$,  $(1,0;s_{R},t_{R})\in\mathcal{B}(T_{v}^{i+1})$, and
 there are seven nonnegative integers $x_{P},x_{S},x_{T},x_{A},y_{S},y_{T},y_{A}$  satisfying  Systems $(1.2)$ or $(1.3)$.
 \end{itemize}
\end{claim}
\begin{proof}
The proofs of the two statements are essentially the same, so we only prove the first one here.  The necessity is obvious.
    Now, suppose  $\phi_{1}$ is a
$(p,q-1;s_{L},t_{L})$-coloring  of $\tilde{T}_{v}^{i}$  and $\phi_{2}$  is
a  $(0,1;s_{R},t_{R})$-coloring of $T_{v}^{i+1}$.
If there are seven nonnegative integers $x_{P},x_{S},x_{T},x_{A},y_{S},y_{T},y_{A}$  satisfying System $(1.1)$, then by exchanging the colors  of $\phi_{2}$ on $T_{v}^{i+1}$ (if necessary) and then combining it with  $\phi_{1}$, one can obtain a satisfactory $(p,q;s,t)$-coloring  of $\tilde{T}_{v}^{i+1}$.
\end{proof}

Before proceeding, we observe the following easy but useful fact, which shows that it suffices to consider only the first seven equalities and inequalities for systems  $(1.1)$, $(1,2)$, $(1.3)$.
\begin{observation}
\label{obser:system}
 In each of Systems $(1.1)$, $(1,2)$, $(1.3)$,   if the $2$nd, $3$rd and $4$th equalities hold, then the $8$th inequality also  holds.
\end{observation}
\begin{proof}
      Assume that the $2$nd, $3$rd and $4$th equalities in System (1.1) hold.
    By summing the $2$nd and $3$rd  equations, we obtain 
$x_A+y_A=s-s_L+t_R-y_T$.
Then, according to the $4$th equation, there is $x_A+y_A=s-s_L+t-t_L$.
Recall that $p+q+s+t\le \Delta$
and 
$a_L=\Delta-p-(q-1)-s_{L}-t_{L}$, we have
$1+x_A+y_A= 1+s-s_L+t-t_L\le 1+\Delta-p-q-s_L-t_L=a_L$ and thus the $8$th inequality holds.
 Similarly,  we can prove this statement also holds for  Systems $(1.2)$ and $(1.3)$. 
\end{proof}

The following claim further reduce  the problem of determining whether $\tilde{T}_{v}^{i+1}$ has a $(p,q;s,t)$-coloring  to
 checking whether the given values satisfy 
the particular set of inequalities (that is, conditions $(1.1)^*$, $(1,2)^*$, $(1.3)^*$ in Lemma~\ref{lemma:iff}). 

\begin{claim}
\label{claim:iff2}
 Given  eight nonnegative integers $p,q,s,t,s_{L},t_{L},s_{R},t_{R}$. Then, there are seven nonnegative integers $x_{P},x_{S},x_{T},x_{A},y_{S},y_{T},y_{A}$  satisfying System $(1.1)$  (resp. $(1.2)$ and $(1.3)$) if and only if  $p,s,t,s_{L},t_{L},s_{R},t_{R}$ satisfy the three conditions in $(1.1)^*$  (resp. $(1.2)^*$ and $(1.3)^*$).
\end{claim}
\begin{proof}
    Suppose first that $x_{P},x_{S},x_{T},x_{A},y_{S},y_{T},y_{A}$ are seven nonnegative integers satisfying System $(1.1)$. The first inequality $t_R\le t \le t_L+t_R$ in $(1.1)^*$ follows directly from the $4$th equality  (implying $ y_T=t_L+t_R-t$) and the $7$th inequality (implying  $0\le y_T\le t_L$) in      System  $(1.1)$.

We prove the second inequality $\max\{0,s_L-s\}\le \min\{s_L,t-t_L\}$ in $(1.1)^*$  by giving the lower and upper bounds of $y_S$.
    By    the $3$rd equality   in   System  $(1.1)$, $y_S=x_A+s_L-s\ge s_L-s$. 
    Since $y_S$ is nonnegative,  we have $y_S\ge \max\{0,s_L-s\}$. 
    Then by the $2$nd equality we have $y_A= t_{R}- y_{S} - y_{T}\ge 0$ and thus  $y_S\le t_R-y_T=t-t_L$ by the $4$th equality.  
    This, together with the $6$th inequality (implying $y_S\le s_L$), implies  that $y_S\le \min\{s_L,t-t_L\}$. Thus, the  second  inequality  in $(1.1)^*$  is true.

Finally, we prove the third inequality $0\le s_L+s_R-s\le s_L+p+t-t_R$   in $(1.1)^*$ .
According to  the $1$st and $3$rd equalities in 
System  $(1.1)$, we have $$s_L+s_R-s=x_{P}+x_{S}+y_{S}+x_{T}\ge0.$$
Then,  by  the $5$th, $6$th and $7$th inequalities, and further by  the $4$th equality  in 
System  $(1.1)$, 
there are
  \[
\begin{array}{rll}
s_L+s_R-s&=x_{P}+x_{S}+x_{T}+y_{S} &  \\
 & \le p+s_L+(t_L-y_T) & \\
 &=s_L+p+t-t_R.
\end{array}
\]  
Thus, the third inequality in $(1.1)^*$ holds.

Conversely, assume that  the seven nonnegative integers $p,s,t,s_{L},t_{L},s_{R},t_{R}$ satisfy the three conditions in  $(1.1)^*$, we define
\begin{enumerate}[label=\ding{\numexpr171+\arabic*}, itemsep=1pt]
    \item  $ y_T=t_L+t_R-t$,
    \item  $y_S=\min\{s_L,t-t_L,s_L+s_R-s\}$,
    \item  $x_A=s-s_L+y_S$, \ (hence $x_A=\min\{s,s-s_L+t-t_L,s_R\}$ by \ding{173}.)
    \item   $y_A=t-t_L-y_S$,  
    \item   $x_T=\min\{t-t_R,s_L+s_R-s-y_S\}$,
    \item  $x_P=\min\{p,s_R-x_A-x_T\}$,
    \ (hence $x_P=\min\{p,s_L+s_R-s-y_S-x_T\}$ by \ding{174}.)
    \item  $x_S=s_R-x_A-x_T-x_P$.
\end{enumerate}
Next, we prove that these seven integers $x_{P},x_{S},x_{T},x_{A},y_{S},y_{T},y_{A}$ satisfy System $(1.1)$.
We first prove that  they are nonnegative one by one as follows.
\begin{itemize}
    \item Due to the first inequality ($t \le t_L+t_R$) in $(1.1)^*$, $y_T\ge0$   is easy to see.
    \item 
Note that $s_L$ is nonnegative. This, together with 
    the conditions $(2)$ ($t-t_L\ge 0$) and $(3)$ ($ s_L+s_R-s\ge0$)
    in  $(1.1)^*$, implies that
 $y_S\ge0$.
\item According to  the condition (2) ($s_L-s\le t-t_L$) in $(1.1)^*$ and $s,s_R\ge0$,  we have $x_A\ge0$.
\item It follows from
\ding{173} that  $y_S\le t-t_L$. Thus  $y_A\ge0$.
\item According to \ding{173}, there is $s_L+s_R-s-y_S\ge0$. This, together with 
the condition $(1)$ ($t\ge t_R$) in   $(1.1)^*$, implies that $x_T\ge0$.
\item 
By \ding{176},  we have $x_T\le s_L+s_R-s-y_S$ and thus $x_P\ge0$.
\item 
By \ding{177}, there is $x_P\le s_R-x_A-x_T$  and so  $x_S\ge0$.
\end{itemize}
Then we prove that the eight conditions in System $(1.1)$ hold one by one as follows.
\begin{itemize}
    \item The $1$st equality $s_{R} = x_{P} + x_{S} + x_{T} + x_{A}$ follows directly from \ding{178}.
    \item It follows from \ding{175} and then \ding{172} that
    $y_{S} +  y_{A}=t-t_L=t_R-y_T$. Thus, the $2$nd equality $t_{R} = y_{S} + y_{T} + y_{A}$ in  System $(1.1)$  holds. 
    \item  The $3$rd equality  $s = s_{L} - y_{S} + x_{A}$ is the same as \ding{174}. 
    \item  The $4$th equality $t = t_{R} + t_{L} - y_{T}$ is the same as \ding{172}.
    \item The $5$th inequality $x_{P} \le p$ is easy to see from \ding{177}.
    \item Due to the condition (3) in   $(1.1)^*$,  it holds that $s_R-s\le t-t_R+p$. Then, by \ding{173}, we have  $y_S\le s_L$ and so   $s_R-s=s_L+s_R-s-s_L\le s_L+s_R-s-y_S$. Therefore,  $$s_R-s\le \min\{t-t_R+p, s_L+s_R-s-y_S\}. $$

Note that for any three numbers  $a,b,c\ge0$ it holds that   $\min\{a,b\}+\min\{c,b-\min\{a,b\}\}=\min\{a+c,b\}$.
This euqality, together with
   \ding{176} and \ding{177}, 
implies that (note that $t-t_R\ge0$, $s_L+s_R-s-y_S\ge0$ by \ding{173} and $p\ge0$)
     \[
\begin{array}{rll}
x_T+x_P&=\min\{t-t_R,s_L+s_R-s-y_S\}+\min\{p,s_L+s_R-s-y_S-x_T\} &  \\
 &=\min\{t-t_R+p,s_L+s_R-s-y_S\}\\
 & \ge s_R-s.
\end{array}
\]  
Finally, according to \ding{178} and \ding{174}, we have 
  \[
\begin{array}{rll}
x_S+y_S&= (s_R-x_A-x_T-x_P)+y_S &  \\
 &=s_R-(s-s_L+y_S)-x_T-x_P+y_S\\
 &=s_L+s_R-s-(x_T+x_P)\\
 &\le s_L+s_R-s-(s_R-s)\\
 &=s_L.
\end{array}
\]  
The $6$th  inequality is proved.
 
\item It follows from \ding{176} that $x_T=\min\{t-t_R,s_L+s_R-s-y_S\} \le t-t_R$.
This, together with \ding{172}, implies  $x_T+y_T\le (t-t_R)+(t_L+t_R-t)=t_L$. Hence, the $7$th inequality holds.

\item The  $8$th  inequality  follows from \Cref{obser:system}.
\end{itemize}

Therefore,  these seven integers $x_{P},x_{S},x_{T},x_{A},y_{S},y_{T},y_{A}$ satisfy System $(1.1)$.

Due to \Cref{obser:system}, we only need to consider the first seven equalities and inequalities for systems  $(1.1)$, $(1,2)$, $(1.3)$.
Now, by regarding $p$
 in System $(1.1)$ as $p-1$, we obtain System $(1,2)$. 
 Similarly,
 by replacing  $p$ as $p-1$ and $s_L$ as $s_L-1$ in System $(1.1)$, we derive System $(1.3)$.
 Hence, the relationships between System $(1.2)$ and the three conditions in  $(1.2)^*$,  and  between System $(1.3)$ and the three conditions in  $(1.3)^*$ are easy to see.
 This finishes the proof.
\end{proof}

According to  \Cref{obser:3cases}, \Cref{claim:iff1}  and \Cref{claim:iff2}, 
Lemma~\ref{lemma:iff} is  true. \qed

\section{Discussion}
\label{sec:summary}

Although semistrong edge coloring is not a novel concept, its study remains largely  unexplored. (This challenge is similar to that of the uniquely restricted edge coloring problem.)
A natural problem is to determine the computational complexity and develop algorithmic approaches for semistrong edge coloring on  certain special classes of graphs.
For example, the problem of deciding whether a graph  admits a strong edge coloring  using $k$ colors
 is known to be NP-complete even for  bipartite graphs and $C_{4}$-free graphs \cite{M2002}. An open question is whether a direct complexity proof can be established specifically for semistrong edge coloring on these two classes.
  Moreover, it is worth exploring whether the algorithm for semistrong edge coloring on trees can be extended to  partial $k$-trees.

The edge coloring problem is closely related to the matching problem. 
In particular,
the maximum induced (strong) matching problem has been extensively  studied from a computational perspective, and polynomial-time algorithms have been developed for several classes of graphs. 
For relevant research on  induced matching, see \cite{golumbic2000new,DDL2013,Joos2016,KANJ2011,CST2003,cameron1989induced}, etc.
As a relaxation of induced matching, 
the concept of semistrong matching arose from the following problem: find the maximum integer $k$ such that $H=kP_3$ (the disjoint union of $k$ copies of the path on 3 vertices) is an induced subgraph of the $n$-dimensional cube $Q^n$, with the additional constraint that each $P_3$  component  has an edge in the same direction of $Q^n$. This problem can be reduced to finding the maximum semistrong matching of $Q^{n-1}$. 
For more details, see \cite{GH2005}, in which the maximum semistrong matching of the subset graphs and the Kneser graphs are studied.
Given the above, the study of semistrong matching,
particularly its computational complexity and algorithmic aspects, deserves further attention.\\

\bigskip

\noindent{\bf Acknowledgements}
We are sincerely grateful to the anonymous reviewers for their careful reading and valuable suggestions, which helped us improve this work.\\

\noindent{\bf Funding}  The first author is supported by the China Scholarship Council (CSC)	and  SEU Innovation Capability Enhancement Plan for  Doctoral Students (CXJH\_SEU 24119).
The second author is supported by NSFC 11771080.\\

\noindent{\bf  Data availability}
No data was used for the research described in the article.\\


\noindent{\bf Declaration of competing interests} The authors declare that they have no conflict of interest.

\bibliographystyle{amsplain}
\bibliography{ref}

@article{HL2017,
  title={On $(s, t)$-relaxed strong edge-coloring of graphs},
  author={He, D. and Lin, W.},
  journal={Journal of Combinatorial Optimization},
  volume={33},
  number={2},
  pages={609--625},
  year={2017},
  publisher={Springer},
  doi={10.1007/s10878-016-0054-9}
}

@incollection{EN1989,
  author    = {Erd\H{o}s, P. and  Ne\v{s}et\v{r}il, J.},
  title     = {Problems},
  booktitle = {Irregularities of Partitions},
  editor    = {G. Hal\'asz and V. T. S\'os},
  pages     = {162--163},
  publisher = {Springer},
  address   = {Berlin},
  year      = {1989}
}

@article{AZ2002,
  title={Algorithmic aspects of acyclic edge colorings},
  author={Alon, N. and Zaks, A.},
  journal={Algorithmica},
  volume={32},
  number={4},
  pages={611--614},
  year={2002},
  publisher={Springer}
}

@phdthesis{G1967,
  author = {Gupta, R.G.},
  title = {Studies in the theory of graphs},
  school = {Tata Institute of Fundamental Research, Bombay},
  year = {1967}
}

@article{GH2005,
  author = {Gyárfás, A. and Hubenko, A.},
  title = {Semistrong edge coloring of graphs},
  journal = {Journal of Graph Theory},
  volume = {49},
  year = {2005},
  pages = {39--47}
}

@article{H1981,
  author = {Holyer, I.},
  title = {The {NP}-completeness of edge colorings},
  journal = {SIAM Journal on Computing},
  volume = {10},
  year = {1981},
  pages = {718--720}
}

@article{LG1983,
  author = {Leven, D. and Galil, Z.},
  title = {{NP} completeness of finding the chromatic index of regular graphs},
  journal = {Journal of Algorithms},
  volume = {4},
  year = {1983},
  pages = {35--44}
}

@article{M2002,
  author = {Mahdian, M.},
  title = {On the computational complexity of strong edge-coloring},
  journal = {Discrete Applied Mathematics},
  volume = {118},
  year = {2002},
  pages = {239--248}
}

@article{V1964,
  author = {Vizing, V.G.},
  title = {On an estimate of the chromatic class of a $p$-graph},
  journal = {Diskretnyi Analiz},
  volume = {3},
  year = {1964},
  pages = {25--30}
}

@article{ASZ2001,
  title={Acyclic edge colorings of graphs},
  author={Alon, N. and Sudakov, B. and Zaks, A.},
  journal={Journal of Graph Theory},
  volume={37},
  number={3},
  pages={157--167},
  year={2001},
  publisher={Wiley Online Library}
}

@article{F1978,
  title={Atsiklicheskij khromaticheskij klass grafa},
  author={Fiam{\v{c}}ik, J.},
  journal={Math. Slovaca},
  volume={28},
  pages={139--145},
  year={1978}
}

@article{LMS2024,
  title={Revisiting semistrong edge-coloring of graphs},
  author={Lu{\v{z}}ar, B. and Mockov{\v{c}}iakov{\'a}, M. and Sot{\'a}k, R.},
  journal={Journal of Graph Theory},
  volume={105},
  number={4},
  pages={612--632},
  year={2024},
  publisher={Wiley Online Library}
}

@article{GHL2001,
  title={Uniquely restricted matchings},
  author={Golumbic, M.C. and Hirst, T. and Lewenstein, M.},
  journal={Algorithmica},
  volume={31},
  pages={139--154},
  year={2001},
  publisher={Springer}
}

@inproceedings{BRS2017,
  title={Uniquely restricted matchings and edge colorings},
  author={Baste, J. and Rautenbach, D. and Sau, I.},
  booktitle={Graph-Theoretic Concepts in Computer Science: 43rd International Workshop, WG 2017, Eindhoven, The Netherlands, June 21-23, 2017, Revised Selected Papers 43},
  pages={100--112},
  year={2017},
  organization={Springer}
}

@article{GHHL2005,
  title={Generalized subgraph-restricted matchings in graphs},
  author={Goddard, W. and Hedetniemi, S.M. and Hedetniemi, S.T. and Laskar, R.},
  journal={Discrete Mathematics},
  volume={293},
  number={1-3},
  pages={129--138},
  year={2005},
  publisher={Elsevier}
}

@article{FJ1983,
  title={Strong edge-colorings of graphs and applications to multi-k-gons},
  author={Fouquet, J.L. and Jolivet, J.L.},
  journal={Ars Combinatoria A},
  volume={16},
  pages={141--150},
  year={1983}
}

@article{golumbic2000new,
  title={New results on induced matchings},
  author={Golumbic, M.C. and Lewenstein, M.},
  journal={Discrete Applied Mathematics},
  volume={101},
  number={1-3},
  pages={157--165},
  year={2000},
  publisher={Elsevier}
}

@article{DDL2013,
title = {New results on maximum induced matchings in bipartite graphs and beyond},
journal = {Theoretical Computer Science},
volume = {478},
pages = {33-40},
year = {2013},
issn = {0304-3975},
author = {Dabrowski, K.k. and Demange, M. and  Lozin, V.V.},
}

@article{Joos2016,
author = {Joos, F.},
title = {Induced Matchings in Graphs of Bounded Maximum Degree},
journal = {SIAM Journal on Discrete Mathematics},
volume = {30},
number = {3},
pages = {1876-1882},
year = {2016}
}

@article{KANJ2011,
title = {On the induced matching problem},
journal = {Journal of Computer and System Sciences},
volume = {77},
number = {6},
pages = {1058-1070},
year = {2011},
author = {Kanj, I. and Pelsmajer, M.J. and  Schaefer, M. and  Xia, G.}
}

@article{CST2003,
title = {Finding a maximum induced matching in weakly chordal graphs},
journal = {Discrete Mathematics},
volume = {266},
number = {1},
pages = {133-142},
year = {2003},
author = {Cameron, K. and Sritharan, R. and Tang, Y.},
}

@article{cameron1989induced,
  title={Induced matchings},
  author={Cameron, K.},
  journal={Discrete Applied Mathematics},
  volume={24},
  number={1-3},
  pages={97--102},
  year={1989},
  publisher={Elsevier}
}

\end{document}